\documentclass[12pt]{amsart}
\usepackage{latexsym, amsmath, amssymb, amsthm}
\usepackage{epsfig}
\usepackage{amscd}
\usepackage{latexsym}
\usepackage{pstricks,pst-node,cite}

\newtheorem{thm}{Theorem}[section]
\newtheorem{lem}[thm]{Lemma}
\newtheorem{prop}[thm]{Proposition}
\newtheorem{cor}[thm]{Corollary}

\theoremstyle{remark}

\newcounter{fignum}
\psset{linecolor=black}\newgray{lightgray}{.9}\newgray{darkgray}{.75}\newgray{pup}{.4}
\psset{linewidth=.6pt,dash=3pt 3pt,doublesep=.05,dotsize=1pt 5}
\SpecialCoor

\def \Aut {\mathrm{Aut}}

\newfont{\sBbb}{msbm7 scaled\magstephalf}

\def\Dist{\hbox{\rm Dist}}
\def\Det{\hbox{\rm Det}}

\begin{document}

\title{The distinguishing numbers of merged Johnson graphs}
\author{Dongseok Kim}
\address{Department of Mathematics \\Kyonggi University
\\ Suwon, 443-760 Korea}
\email{dongseok@kgu.ac.kr}

\author{Young Soo Kwon}
\address{Department of Mathematics \\Yeungnam University \\Kyongsan, 712-749, Korea}
\email{yskwon@yu.ac.kr}
\thanks{This research was supported by the Yeungnam University
research grants in 2010.}

\author{Jaeun Lee}
\address{Department of Mathematics \\Yeungnam University \\Kyongsan, 712-749, Korea}
\email{julee@yu.ac.kr}

\date{}
\maketitle

\begin{abstract}
In present article, we determine
the distinguishing number of the merged Johnson graphs which are
generalization of both the Kneser graphs and of the Johnson graphs. 
\end{abstract}

\section{Introduction}

The \emph{distinguishing number} of a graph $G$ is the minimum
number of colors for which there exists an assignment of colors to
the vertices of $G$ such that the identity is the only
color-preserving automorphism of $G$.  Generally, for a
permutation group $\Gamma$ acting on $X$, the \emph{distinguishing
number} of $\Gamma$ is the minimum number of cells of a partition
$\pi$ of $X$ satisfying that the identity is the only element of
$\Gamma$ fixing each cell of $\pi$. Albertson and
Collins first introduced distinguishing number of a graph~\cite{AC} and there have
been many interesting results on the distinguishing numbers of graphs and
permutation groups in last few years~\cite{AB, AB1, AC, BC, C, FI, IK, KZ, ZW}.

Here we consider a class of graphs based on the Johnson graphs.
For positive integers $k, n$ such that $1 \le k \le \frac{n}{2}$,
the \emph{Johnson graph $J(n,k)$} has vertices given by the $k$-subsets
of $[n]=\{1,2,\ldots,n\}$ and there exists an edge between two vertices if and only
if their intersection has size $k-1$. Given a nonempty subset $I
\subseteq \{1,2,\ldots, k\}$, the \emph{merged Johnson graph $J(n,k)_I$}
is the union of the distance $i$ graphs $J(n,k)_i$ of $J(n,k)$ for
all $i$, namely, two $k$-subsets are adjacent in $J(n,k)_I$ if and
only if their intersection has $k-i$ elements for some $i \in I$.
The merged Johnson graphs $J(n,k)_I$ include many interesting graphs such as the
Johnson graph $J(n,k)=J(n,k)_{\{1\}}$ and the Kneser graph
$K(n,k)=J(n,k)_{\{k\}}$.

In \cite{AB}, M.O. Albertson and D.L. Boutin determined the
distinguishing number of the Kneser graphs. The aim of the present article is to
determine the distinguishing number of the merged Johnson graphs.

The outline of this paper is as follows. In section~\ref{pre},
we review some preliminaries regarding the distinguishing numbers and
the merged Johnson graphs.
In section~\ref{mainsec}, we find Theorem~\ref{main} which addresses the distinguishing numbers of the merged Johnson graphs.
We also prove lemmas which are used for a proof of the main theorem.
At last, we provide a proof of Theorem~\ref{main} in section~\ref{proof}.

\section{Preliminaries} \label{pre}

For a given graph $G$, a coloring $f: V(G) \rightarrow \{1,2,\ldots,r\}$ is said to be  \emph{$r$-distinguishing}
if the identity is the only graph automorphism
$\phi$ satisfying $f(\phi(v))=f(v)$ for all $v \in V(G)$. This
means that the distinguishing coloring is a symmetry-breaking
coloring of $G$. The distinguishing number, denoted by $\Dist(G)$,
is the minimum $r$ that $G$ has an $r$-distinguishing coloring.
One can easily see that $\Dist(G)=\Dist(G^c)$ where $G^c$ is the complement of
$G$. If $G$ is an \emph{asymmetric graph}, the identity is the
only automorphism of $G$, then $\Dist(G)=1$. In fact, the converse
is also true.

For a graph $G$ and for a subset  $S\subseteq V(G)$, a coloring
$f: S \rightarrow \{1, 2, \ldots, s\}$ is said to be
\emph{$s$-distinguishing} if for any graph automorphism $\phi \in
\Aut(G)$ fixing $S$ set-wisely and if $f(\phi(v))=f(v)$ for all $v\in
S$ then $\phi$ fixes all elements of $S$. In this case, $\phi$ does
not need to fix other vertices outside of the given set $S$. If there exists an
$s$-distinguishing coloring for $S$, the set $S$ is called an
\emph{$s$-distinguishable set}.

For a graph $G$, a subset $S\subseteq V(G)$ is called a \emph{determining set}
if the identity is the only automorphism fixing every element of $S$.
The \emph{determining number} of $G$, denoted by $\Det(G)$, is the
minimum number $t$ that $G$ has a determining set of the cardinality
$t$. It is also true that $\Det(G)=\Det(G^c)$. For example, for any
$n \ge 3$, $\Det(C_n)=2$ because the set of two adjacent vertices
is a determining set of $C_n$ and there is no determining set of
$C_n$ of the cardinality $1$. Note that if $S\subseteq V(G)$ is a
determining set, then any subset $T \subseteq V(G)$ containing $S$ is also a determining set.     
The determining sets provide a useful tool for finding the distinguishing number of $G$ as
stated in the following theorem.

\begin{thm} {\rm (\cite{AB})}  \label{dist-det}
For a given graph $G$, $G$ has an $r$-distinguishable determining set if and only if $G$ has a $(r+1)$-distinguishing coloring.
\end{thm}

Consequently, we find the following corollary which will be used in the proof of our results.

\begin{cor}   \label{asymm-dist}
 For a given graph $G$, if there is a determining set $S$ of $G$ such that the induced subgraph of $S$ is asymmetric then
 $\Dist(G)=1$ or $2$.
\begin{proof}
For a graph automorphism $\phi$ of $G$ fixing $S$ set-wisely, the
restriction of $\phi$ on the induced subgraph $\langle S \rangle$
is a graph automorphism of $\langle S \rangle$. Since $\langle S
\rangle$ is asymmetric, the coloring $f(v)=1$ for all $v \in S$ is
a $1$-distinguishing. Furthermore, since $S$ is a determining set,
$G$ has $2$-distinguishing coloring by Theorem \ref{dist-det}.
\end{proof}
\end{cor}

Using Theorem \ref{dist-det}, M.O. Albertson and D.L. Boutin determined the distinguishing number of the Kneser graphs as follows.

\begin{prop} {\rm (\cite{AB})}  \label{dist-Kneser}
 For any integers $n > k \ge 2$ with $n > 2k$, $\Dist(K(n,k))=2$ except $(n,k)=(5,2)$.
\end{prop}

Note that the Kneser graph $K(5,2)$ is isomorphic to the Pertersen graph and its distinguishing number is $3$.
$K(n,1)$ is isomorphic to the complete graph and so $\Dist(K(n,1))=n$ for any integer $n$.

For a graph $G$, the distinguishing number $\Dist(G)$ equal to the distinguishing number of $\Aut(G)$ which acts on
the vertex set $V(G)$. Hence we have the following lemma. The proof is straightforward and we omit it.

\begin{lem} \label{same-dist}
Let  $G_1$ and $G_2$ be  two graphs having the same vertex set
$V$.
\begin{enumerate}
\item [{\rm (1)}] If $\Aut(G_1)$ is a subgroup of $\Aut(G_2)$ as
acting groups on $V$, then $\Dist(G_1) \le \Dist(G_2)$. \item
[{\rm (2)}] If $\Aut(G_1)=\Aut(G_2)$ as acting groups on $V$, then
$\Dist(G_1) = \Dist(G_2)$.
\end{enumerate}
\end{lem}

\section{The Distinguishing numbers of the merged Johnson graphs} \label{mainsec}

Let $\Omega$ be the set of all $k$-subsets of $[n]$. The action of $S_n$ on $[n]$ naturally induces an action of $S_n$ on $\Omega$.
Then the Johnson graph $J(n,k)$ is an orbital graph which corresponds to the orbital $\{ (M,N)\in \Omega^2 \ | \ |M \cap N|=k-1 \}$,
that is, the orbit of $S_n$ on $\Omega^2$. The \emph{distance $i$ graph $J(n,k)_i$} of
$J(n,k)$ is an orbital graph corresponding to the orbital
$$\Gamma_i = \{ (M,N)\in \Omega^2 \ | \ |M \cap N|=k-i \}$$
(see \cite{DM} for orbital graphs).

For a merged Johnson graph $J=J(n,k)_I$, the complementation of each $k$-sets $M \rightarrow M^c$ induces an isomorphism from
$J(n,k)_I$ to $J(n,n-k)_I$. So we may assume  that $k \le n/2$.
For a merged Johnson graph $J=J(n,k)_I$ with $1 \le k \le \frac{n}{2}$, if $I=\emptyset$ or $\{1$, $2$, $\ldots$, $k\}$ then
$J$ is the null or complete graph and so $\Aut(J)=S_d$, where $d= {n \choose k}$. Thus, we further assume that $k \ge 2$ and
$\emptyset \varsubsetneq I \varsubsetneq \{1$, $2$, $\ldots$, $k\}$. For
notational simplicity, let $I' = I \setminus \{k\}$, and for any
integer $t$, let $t-I = \{t-i \ | \ i \in I \}$ and $t-I' = \{t-i
\ | \ i \in I' \}$. We also denote $I'' = k - I'$, and let
$e=\frac{1}{2} {n \choose {n/2}}$.

In \cite{J}, G. Jones found the automorphism groups of the merged Johnson graphs as follows.

\begin{thm} {\rm (\cite{J})}  \label{Aut-merged}
Let $J=J(n,k)_I$, where  $2 \le k \le \frac{n}{2}$ and
$\emptyset \varsubsetneq I \varsubsetneq \{1,2,\ldots,k\}$ and let $A=\Aut(J)$.
\begin{enumerate}
\item [{\rm (1)}] If $2 \le k < \frac{n-1}{2}$, and $J \neq
J(12,4)_I$ with $I= \{1,3 \}$ or  $\{2,4 \}$, then $A=S_n$ with
orbitals $\Gamma_0, \Gamma_1, \ldots, \Gamma_k \subset \Omega^2$.
\item [{\rm (2)}] If $J = J(12,4)_I$ with $I= \{1,3 \}$ or  $\{2,4
\}$, then $A=O^{-1}_{10}(2)$  with orbitals $\Gamma_0, \Gamma_1
\cup \Gamma_3, \Gamma_2 \cup \Gamma_4$. \item [{\rm (3)}] If
$k=\frac{n-1}{2}$ and $I \neq k+1-I$, then $A=S_n$ with orbitals
$\Gamma_0, \Gamma_1, \ldots, \Gamma_k \subset \Omega^2$. \item
[{\rm (4)}] If $k=\frac{n-1}{2}$ and $I = k+1-I$, then $A=S_{n+1}$
with orbitals $\Gamma_0$ and $\Gamma_i \cup \Gamma_{k+1-i}$ for
all $i=1,2, \ldots, \lfloor \frac{k+1}{2} \rfloor$. \item [{\rm
(5)}] If $k=\frac{n}{2}$ and $I \neq \{k\}$ nor
$\{1,2,\ldots,k-1\}$, and $I' \neq I''$, then $A=S_2 \times S_n$
with orbitals $\Gamma_0, \Gamma_1, \ldots, \Gamma_k \subset
\Omega^2$. \item [{\rm (6)}] If $k=\frac{n}{2}$ and $I \neq \{k\}$
nor $\{1,2,\ldots,k-1\}$, and $I' = I''$, then $A=S^{e}_{2}:S_n$
with orbitals $\Gamma_0$ and $\Gamma_i \cup \Gamma_{k-i}$ for all
$i=1,2, \ldots, \lfloor \frac{k}{2} \rfloor$ and $\Gamma_k$. \item
[{\rm (7)}] If $k=\frac{n}{2}$ and $I = \{k\}$ or
$\{1,2,\ldots,k-1\}$, then $A=S_{2}^{e}:S_e = S_2 \wr S_e$ with
orbitals $\Gamma_0, \Gamma_1\cup \cdots \cup \Gamma_{k-1}$ and
$\Gamma_k$.
\end{enumerate}
\end{thm}

To understand the automorphism group $S_{n+1}$ of $J(n,
\frac{n-1}{2})_I$ with $I = k+1-I$, let $\tilde{[n]} = [n] \cup \{
\infty \}$ and let $\Psi$ be the set of equipartitions of
$\tilde{[n]}$, by which we mean the unordered partitions
$\{P_1,P_2\}$ of $\tilde{[n]}$ satisfying $|P_1|=|P_2|=
\frac{n+1}{2}$. There is a bijection $\phi: \Omega \rightarrow
\Psi$, sending each $M$ to $\{ M \cup \{ \infty \} , [n]-M \}$.
Note that its inverse sends an equipartition $\{P_1,P_2\}$ to $P_i
\setminus \{ \infty \}$, where $i$ is chosen so that $\infty \in
P_i$. The natural action of $S_{n+1}$ on $\tilde{[n]}$ induces an
action of $S_{n+1}$ on $\Omega$. By the condition $I = k+1-I$, one
can see that this action induces an automorphism group of $J(n,
\frac{n-1}{2})_I$ (For a detail information, see the paper\cite{J}).
The next theorem is the main theorem of this paper.

\begin{thm}  \label{main}
Let $J=J(n,k)_I$, where  $1 \le k \le \frac{n}{2}$ and
$\emptyset \varsubsetneq I \varsubsetneq \{1,2,\ldots,k\}$.
\begin{enumerate}
\item [{\rm (1)}] If $k=1$ then  $\Dist(J)= n$. \item [{\rm (2)}]
If $(n,k) \neq (5,2)$ and $J\neq J(n, \frac{n}{2})_I$  with $I =
\{ \frac{n}{2} \}$, $\{1,2,\ldots,\frac{n}{2}-1\}$ or $I' = I''$,
then $\Dist(J) = 2$. \item [{\rm (3)}] If $J = J(5,2)_I$ with $I =
{\{1\}}$ or $\{ 2 \}$; or  $J= J(n, \frac{n}{2})_I$ satisfying $I'
= I''$ and $I$ is neither $\{\frac{n}{2}\}$ nor
$\{1,2,\ldots,\frac{n}{2}-1\}$, then $\Dist(J)=3$. \item [{\rm
(4)}] If $J= J(n, \frac{n}{2})_I$ with $I = \{\frac{n}{2}\}$ or
$\{1,2,\ldots,\frac{n}{2}-1\}$, then $\Dist(J) = \lceil
\frac{1+\sqrt{1+4 {n \choose n/2}}}{2} \rceil$.
\end{enumerate}
\end{thm}

\begin{cor}  \label{Dist-Johnson}
Let $G=J(n,k)$ be the Johnson graph. Then
\begin{enumerate}
\item [{\rm (1)}] if $k=1$, then  $\Dist(G)= n$; \item [{\rm (2)}]
if $(n,k) \neq (4,2)$ nor $(5,2)$ , then $\Dist(G) = 2$ and \item
[{\rm (3)}] if $G = J(4,2)$ or $J(5,2)$, then $\Dist(J)=3$.
\end{enumerate}
\end{cor}

We will prove Theorem \ref{main} in the next section. For the rest of the section, we will prove lemmas which will be
used in the proof of Theorem \ref{main}. For a set
$X=\{ 1$, $2$, $3$, $\ldots$, $n \}$ and for any permutation $\phi$ of $X$,
$\phi$ can be represented by $(i_1$, $i_2$, $\ldots$, $i_n)$, where
for any $j =1$, $2$, $\ldots$, $n$, $\phi(j) = i_j$. Throughout the rest of the paper, we use
the above representation of permutations. For any permutation
$\phi = (i_1$, $i_2$, $\ldots$, $i_n)$ and for any $k, \ell \in [n]$
with $1 \le k \le \frac{n}{2}$, let $V^{\phi}_{\ell}$ be the
$k$-subset $\{ i_\ell$, $i_{\ell+1}$, $\ldots$, $i_{\ell+k-1} \}$,
where the subscripts are considered as their residue classes
modulo $n$. For our convenience, if $\phi$ is the identity, we
denote $V^{\phi}_{\ell}$ simply by $V_{\ell}$. For any
$i=1,2,\ldots,n$, let $\tau_i$ be the transposition of $X$
exchanging $i$ and $i+1$.

\begin{lem} \label{12-4-exchange}
Let $J=J(12,4)_I$ with $I= \{1,3 \}$. For any permutation $\phi$
of $X$, if an automorphism $\Psi$ of $J$ fixes all vertices in $S
= \{ V^{\phi}_{j}, V^{\phi\tau_1}_{j}, V^{\phi\tau_2}_{j} \ | \
j=1,2,\ldots, 12 \}$, then $\Psi$ also fixes $V^{\phi\tau_i}_j$
for all $i,j=1,2,\ldots,12$.
\begin{proof}
Without any loss of generality, we may assume that $\phi$ is the
identity. Note that $S = \{V_j \ | \ j=1$, $2$, $\ldots$, $12\}$ $\cup$
$\{V^{\tau_1}_2$, $V^{\tau_1}_{10}\}$ $\cup$ $\{V^{\tau_2}_3$, $
V^{\tau_2}_{11}\}$, where $V^{\tau_1}_2 = \{1$, $3$, $4$, $5\},
V^{\tau_1}_{10}$ $= \{10$, $11$, $12$, $2 \}$ and $V^{\tau_2}_3$ $= \{2$, $4$, $5$, $6
\}$, $V^{\tau_2}_{11}$ $= \{11$, $12$, $1$, $3 \}$.

Let $\Psi$ be an automorphism of $J$ fixing all vertices in $S$.
Note that $\{ V^{\tau_3}_{j} \ | \ j=1$, $2$, $\ldots$, $12 \} - S$ $= \{
V^{\tau_3}_{4}$, $V^{\tau_3}_{12} \}$ because $\tau_3$ $=
(1$, $2$, $4$, $3$, $5$, $6$, $\ldots$, $12)$, where $V^{\tau_3}_{4}$ $= \{3$, $5$, $6$, $7 \}$
and $V^{\tau_3}_{12}$ $= \{12$, $1$, $2$, $4 \}$. Let $A$ and $B$ be the sets
of all vertices in $S$ which are adjacent to $V^{\tau_3}_{4}$ and
$V^{\tau_3}_{12}$, respectively. Then we have
$$A = \{ V_1, V_3, V_4, V_5, V_{7}, V_{12}, V^{\tau_2}_{11}\} \ \
\mbox{and} \ \ B = \{ V_1 , V_3, V_4, V_9, V_{11}, V_{12} \}.$$

Let $V$ be a vertex whose all adjacent vertices in $S$
is $A$. For the first case, assume that $|V \cap V_1| = 3$. If $V
\cap V_1 = \{1,2,3 \}$, then $4$, $5$, $6 \notin V$ and $7 \in V$ because
$V$ is adjacent to both $V_3$ and $V_4$ but not to $V_2$. In this
case, $V$ is adjacent to $V_6$, a contradiction. If  $V \cap V_1 =
\{1$, $2$, $4 \}$, then $5$, $6$, $7 \notin V$ and $8 \in V$, and hence $V$ is
adjacent to $V_6$, a contradiction. Similarly, one can show that
$V \cap V_1$ is neither $\{1$, $3$, $4 \}$ nor $\{2$, $3$, $4 \}$. Therefore,
we have $|V \cap V_1| = 1$.

By considering the fact that $V$ is adjacent to $V_1$, $V_3$, $V_4$,
$V_5$, $V_{7}$, $V_{12}$ but not to $V_2$, $V_6$, $V_8$, $V_9$, $V_{10}$, $
V_{11}$, one can show that $V = \{1$, $6$, $9$, $11\}, \{2$, $5$, $10$, $11\}$ or
$\{3$, $5$, $6$, $7\}$. Since $V$ is not adjacent to $V^{\tau_1}_2 =
\{1$, $3$, $4$, $5 \}$, $V$ is $V^{\tau_3}_{4} = \{3$, $5$, $6$, $7\}$.  This
implies that $\Psi$ also fixes $V^{\tau_3}_{4}$.

Let $V'$  be a vertex whose all adjacent vertices in
$S$ is $B$. By considering the fact that $V'$ is adjacent to $V_1$
, $V_3$, $V_4$, $V_9$, $V_{11}$, $V_{12}$ but not to $V_2$, $V_5$,
$V_6$, $V_7$, $V_8$, $V_{10}$, one can show that $V = \{1$, $6$, $8$, $10\},
\{2$, $5$, $8$, $9\}$ or $\{12$, $1$, $2$, $4\}$. Since $V'$ is not adjacent to
$V^{\tau_1}_2 = \{1$, $3$, $4$, $5 \}$, $V'$ is $V^{\tau_3}_{12} =
\{12$, $1$, $2$, $4\}$.  This implies that $\Psi$  fixes $V^{\tau_3}_{12}$.
Up to now, we showed that $\Psi$  fixes $V^{\tau_3}_j$ for all
$j=1$, $2$, $\ldots$, $12$.

Since $\Psi$ fixes all elements in $\{ V_{j}$, $V^{\tau_2}_{j}$,
$V^{\tau_3}_{j} \ | \ j=1$, $2$, $\ldots$, $12 \}$, one can show that
$\Psi$  fixes $V^{\tau_4}_j$ for all $j=1$, $2$, $\ldots$, $12$ by a similar
way. Continuing the similar process, one can show that $\Psi$
fixes $V^{\tau_i}_j$ for all $i$, $j=1$, $2$, $\ldots$, $12$.
\end{proof}
\end{lem}

\begin{lem} \label{12-4-transposition}
Let $J=J(12,4)_I$ with $I= \{1,3 \}$. For any permutation $\phi$
of $X$, if an automorphism $\Psi$ of $J$ fixes all vertices in
$S_1 = \{ V^{\phi}_{j}, V^{\phi\tau_i}_{j} \ | \ i,j=1,2,\ldots,
12 \}$, then $\Psi$ also fixes $V^{\phi\tau_i\tau_k}_j$ for all
$i,j,k=1,2,\ldots,12$.
\begin{proof}
Without any loss of generality, assume that $\phi$ is the
identity. Let $\Psi$ be an automorphism of $J$ fixing all vertices
in $S_1$. Note that $\{ V^{\tau_1\tau_2}_{j} \ | \ j=1$, $2$, $\ldots$,
$12 \} - S_1 = \{ V^{\tau_1\tau_2}_{3}$, $V^{\tau_1\tau_2}_{11} \}$
because $\tau_1\tau_2 = (2$, $3$, $1$, $4$, $5$, $\ldots$, $12)$, where
$V^{\tau_1\tau_2}_{3} = \{1$, $4$, $5$, $6 \}$ and $V^{\tau_1\tau_2}_{11} =
\{11$, $12$, $2$, $3 \}$.  Let $C$ and $D$ be the sets of all vertices in
$S_1$ which are adjacent to $V^{\tau_1\tau_2}_{3}$ and
$V^{\tau_1\tau_2}_{11}$, respectively. Then we have
$$\{ V_3, V_4, V_6, V_{10}, V_{11}, V_{12}, V^{\tau_1}_{2}, V^{\tau_2}_{3}\} \subset C \ \
\mbox{and} \ \  \{ V_3 , V_8, V_{11}, V_{12},
V^{\tau_1}_{2}, V^{\tau_2}_{3} \} \subset D.$$

Let $V$ be a vertex whose all adjacent vertices in
$S_1$ is $C$. By considering the fact that $V$ is adjacent to $V_3$, $V_4$, $V_6$, $V_{10}$, $V_{11}$, $V_{12}$ but not to $V_1$, $V_2$, $V_5$,
$V_7$, $V_8$, $V_9$, we have that $V = \{6$, $8$, $9$, $12\}, \{1$, $4$, $5$, $6\}$ or
$\{2$, $4$, $9$, $10\}$. Since $V$ is  adjacent to both $V^{\tau_1}_2 =
\{1$, $3$, $4$, $5 \}$ and $V^{\tau_2}_3 = \{2$, $4$, $5$, $6 \}$, $V$ is
$V^{\tau_1\tau_2}_{3} = \{1$, $4$, $5$, $6\}$.  This implies that $\Psi$
also fixes $V^{\tau_3}_{4}$.

Let $V'$ be a vertex whose all adjacent vertices in
$S_1$ is $B$. By considering the fact that $V'$ is adjacent to
$V_3$, $V_8$, $V_{11}$, $V_{12}$ but not to $V_1$, $V_2$, $V_4$, $V_5$,
$V_6$, $V_7$, $V_9$, $V_{10}$, one can show that $V = \{6$, $7$, $10$, $12\},
\{2$, $4$, $7$, $8\}$ or $\{11$, $12$, $2$, $3\}$. Since $V'$ is adjacent to both
$V^{\tau_1}_2 = \{1$, $3$, $4$, $5 \}$ and $V^{\tau_2}_3 = \{2$, $4$, $5$, $6 \}$,
$V'$ is $V^{\tau_1\tau_2}_{11} = \{11$, $12$, $2$, $3 \}$.  This implies
that $\Psi$  fixes $V^{\tau_1\tau_2}_{11}$. Therefore $\Psi$ fixes
$V^{\tau_1\tau_2}_j$ for all $j=1$, $2$, $\ldots$, $12$.

Since $\Psi$ fixes all vertices in $ \{ V^{\phi}_{j}$,
$V^{\phi\tau_1}_{j}$, $V^{\phi\tau_2}_{j} \ | \ j=1$, $2$, $\ldots$, $12 \}$
with $\phi = \tau_1$, $\Psi$ also fixes $V^{\tau_1\tau_i}_j$ for
all $i$, $j=1$, $2$, $\ldots$, $12$ by Lemma \ref{12-4-exchange}. By a similar way, one can show that $\Psi$  fixes $V^{\tau_i\tau_k}_j$ for all
$i$, $j$, $k=1$, $2$, $\ldots$, $12$.
\end{proof}
\end{lem}

\begin{lem} \label{12-4--det}
Let $J=J(12,4)_I$ with $I= \{1$, $3 \}$. Then $S = \{ V_{j}$,
$V^{\tau_1}_{j}$, $V^{\tau_2}_{j} \ | \ j=1$, $2$, $\ldots$, $12 \}$ is a
determining set.
\begin{proof} Let $\Psi$ be an automorphism of $J$ fixing all
vertices in  $S = \{ V_{j}, V^{\tau_1}_{j}, V^{\tau_2}_{j} \ | \
j=1,2,\ldots, 12 \}$. By Lemma \ref{12-4-exchange}, $\Psi$  fixes
$V^{\tau_i}_j$ for all $i$, $j=1$, $2$, $\ldots$, $12$. Furthermore, $\Psi$
 fixes $V^{\tau_i\tau_k}_j$ for all $i$, $j$, $k=1$, $2$, $\ldots$, $12$ by Lemma
 \ref{12-4-transposition}. By applying Lemma
 \ref{12-4-transposition} again with $\phi = \tau_i$, one can show
 that $\Psi$ fixes $V^{\tau_i\tau_k\tau_{\ell}}_j$ for all $i,j,k, \ell
 =1,2,\ldots,12$. Continuing the similar process, one can show that
 $\Psi$ fixes $V^{\tau_{i_1}\tau_{i_2}\cdots \tau_{i_t}}_j$ for any positive
 integer $t$ and for all
 $1 \le i_1$, $i_2$, $\ldots$, $i_t$, $j \le 12$. Since $\{ \tau_i \ | \ i=1$, $2$, $\ldots$, $12 \}$ generates symmetric group on
 $X=\{ 1$, $2$, $3$, $\ldots$, $12 \}$, $\Psi$ fixes all permutation of $X$, $i.e.$, $S$ is a  determining set.
\end{proof}
\end{lem}

Let $J$ be a merged Johnson graph $J(2m,m)_I$ with $I \subseteq \{1$, $2$, $\ldots$, $m \}$. For any $v \in V(G)$, let $\bar{v}$ be the vertex $[2m]-v$
for our convenience.

\begin{lem} \label{n,k-det}
Let $J=J(n,k)_I$, where  $1 \le k \le \frac{n}{2}$ and $\emptyset
\varsubsetneq I \varsubsetneq \{1$, $2$, $\ldots$, $k\}$.
\begin{enumerate}
\item[{\rm (1)}]  For $(n,k)=(2m+1,m)$ with $m \ge 3$ and  $I=
\{1,m \}$, $S_1 = \{V_{1}$, $V_{2}$, $\ldots$, $V_{m+2} \}$ is a
determining set.
 \item [{\rm (2)}]For $(n$, $k)=(2m$, $m)$ with $m \ge 3$ and  $I= \{1 \}$, $S_2 = \{V_{1}$, $V_{2}$, $\ldots$, $V_{2m} \}$ $\cup$ $\{  \{1$, $2$, $\ldots$, $m-2$, $m$, $m+2 \} \}$ is a determining set.
\end{enumerate}
\begin{proof} (1) For any automorphism $\phi$ of $J$ as a permutation on vertices of $J$,
let $\phi'$ be a corresponding permutation of  $[n] \cup \{ \infty
\}$. Let $\phi$ be an automorphism of $J$ fixing all elements in
$S_1$. Since $\phi$ fixes $V_1$, $\phi'$ fixes $\{1$, $2$, $\ldots$, $m$,
$\infty \}$ set-wisely or $\phi'$ sends $\{1$, $2$, $\ldots$, $m$, $\infty \}$
to $\{ m+1$, $m+2$, $\ldots$, $2m+1 \}$ set-wisely.

\medskip

\noindent Case 1: $\phi'$ fixes $\{1,2,\ldots, m , \infty  \}$
set-wisely.

Since $\phi$ fixes $V_2$ and $\phi'$ fixes $\{1,2,\ldots, m ,
\infty \}$ set-wisely, $\phi'$ also fixes $\{2,3,\ldots, m+1 , \infty
\}$ set-wisely. This implies that $\phi'$ fixes both $1$ and $m+1$.
Using the fact that $\phi$ fixes all elements in $S_1$, one can
see that $\phi'$ fixes all elements in $[\widetilde{2m+1}] =
[2m+1] \cup \{ \infty \}$, namely, $\phi$ is the identity element.

\medskip

\noindent Case 2: $\phi'$ sends $\{1,2,\ldots, m , \infty  \}$ to $\{
m+1 , m+2 , \ldots, 2m+1 \}$ set-wisely.

Since $\phi$ fixes $V_2$ and $\phi'$  sends $\{1,2,\ldots, m ,
\infty  \}$ to $\{ m+1 , m+2 , \ldots, 2m+1 \}$ set-wisely, $\phi'$
also sends $\{2,3,\ldots, m+1 , \infty  \}$ to $\{ m+2 , m+3 ,
\ldots, 2m+1,1 \}$ set-wisely. By the similar way, one can show that
for any $i=1,2,\ldots, m+2$, $\phi'$ sends $\{i,i+1,\ldots,m+i-1 ,
\infty \}$ to $[\tilde{n}] - \{i,i+1,\ldots,m+i-1   \}$ set-wisely.

Let $\phi'(\infty ) = a$. Then $a$ is contained to $\{ m+1 , m+2 ,
\ldots, 2m+1 \}$ and hence $\phi'$ can not send
$\{a-m+1,a-m+2,\ldots,a , \infty  \}$ to  $[\widetilde{2m+1}] -
\{a-m+1,a-m+2,\ldots,a \}$ set-wisely, which is a contradiction.

Therefore $S_1 = \{V_{1}, V_{2}, \ldots, V_{m+2} \}$ is a
determining set.

\medskip

\noindent (2) Let $\alpha$ be the automorphism which sends $v$ to
$\bar{v}$ for all $v \in V(G)$. Then the order of $\alpha$ is 2
and $Aut(J) \cong  \langle \alpha \rangle  \times S_n$. Note that
for any automorphism $\psi$ of $J$, $\psi$ is either an
automorphism induced by a permutation of $[n]$
 or a product of $\alpha$ and an automorphism induced by a permutation of $[n]$.

Let $\psi$ be an automorphism of $J$ fixing all elements in
$S_2$. Since $\psi$ fixes $V_1$, $\psi$ is an automorphism induced by a permutation of $[n]$
fixing $\{1$, $2$, $\ldots$, $m\}$ set-wisely or $\psi$ is a product of $\alpha$ and an automorphism induced by a permutation of $[n]$
sending $\{1$, $2$, $\ldots$, $m\}$ to  $\{ m+1$, $m+2$, $\ldots$, $2m \}$ set-wisely.

If $\psi$ is an automorphism induced by a permutation $\psi'$ of
$[n]$ fixing $\{1$, $2$, $\ldots$, $m  \}$ set-wise, then $\psi'$ also
fixes $\{i$, $i+1$, $\ldots$, $i+m-1 \}$ set-wise for all $i=1$, $2$, $\ldots$,
$2m$ because $\psi$ fixes $V_i$ for all $i=1$, $2$, $\ldots$, $2m$. This
implies that $\psi$ is the identity. Hence we can assume that $\psi$ is a
product of $\alpha$ and an automorphism induced by a permutation
$\psi''$ of $[n]$ sending $\{1$, $2$, $\ldots$, $m \}$ to  $\{ m+1$, $m+2$, $\ldots$, $2m \}$ set-wisely. Since $\psi$ fixes $V_2$, $\psi''$ also
sends $\{2$, $3$, $\ldots$, $m+1 \}$ to  $\{ m+2$, $m+3$, $\ldots$, $2m$, $1
\}$ set-wisely. This implies that $\psi''$ exchanges $1$ and $m+1$.
By a similar way, one can show that $\psi''$ exchanges $i$ and
$m+i$ for any $i=1$, $2$, $\ldots$, $m$. But in this case, $\psi$ does
not fix $\{1$, $2$, $\ldots$, $m-2$, $m$, $m+2 \}$. Therefore,  $S_2$ is a
determining set.
\end{proof}
\end{lem}

\begin{lem} \label{n=2m-m-Dist>2}
For $m \ge 4$, let $J=J(2m,m)_I$ with $I= \{1,m-1 \}$. Then,
$\Dist(J) > 2$.
\begin{proof}  For any $v \in V(J)$, let $\beta_v$ be the automorphism of $V(J)$ exchanging $v$ and $\bar{v}$ and fixing all other vertices.
Let $f: V(J) \rightarrow \{1,2 \}$ be a coloring. If there exists
$u \in V(J)$ such that $f(u) = f( \bar{u})$ then $\beta_u$ is a
color-preserving automorphism, and hence $f$ is not
$2$-distinguishing. Assume that for all $v \in V(J)$, $f(v)$ and
$f(\bar{v})$ are distinct. Let $\phi$ be an automorphism of $J$
induced by a non-identity permutation of $[2m]$. Let
$$\psi =\left( \prod_{ \{u,\bar{u} \},~  f(u)\neq f(\phi(u)) } \beta_u \right) \phi.$$
Then $\psi$ is a color-preserving automorphism. Therefore, there is
no $2$-distinguishing coloring, and hence  $\Dist(J) > 2$.
\end{proof}
\end{lem}

For $n=2m$ with $m \ge 4$, let $\Phi$ be the set of equipartitions
of $[n]$. Note that the size of $\Phi$ is $\frac{1}{2} {n \choose
m}$. The natural action of $S_{n}$ on $[n]$ induces an action of
$S_n$ on $\Phi$. Let $\alpha$ be a permutation satisfying that for
all non-identity permutation  $\beta \in S_n$, the number of
equipartitions fixed by $\alpha$ in  $\Phi$ is greater than or
equal to the number of equipartitions fixed by $\beta$ in  $\Phi$.
Suppose that there exist a permutation $\gamma \in S_n$ and $i_1,
i_2 , \ldots, i_t \in [n]$ with $t \ge 3$  such that
$\gamma(i_j)=i_{j+1}$ for all $j = 1$, $2$, $\ldots$, $t-1$ and
$\gamma(i_t)=i_1$. Let $\tilde{\gamma} \in S_n$ be a permutation
that $\tilde{\gamma}(i_2 ) = i_1$, $\tilde{\gamma}(i_t ) = i_3$
and $\tilde{\gamma}(\ell) = \gamma(\ell)$ for all $\ell \in [n]
\setminus \{i_2$, $i_t\}$. Then, all equipartitions in  $\Phi$
fixed by $\gamma$ are also fixed by $\tilde{\gamma}$. This implies
that $\alpha$ is a product of disjoint transpositions. For two
permutations $\gamma_1$, $\gamma_2$ $\in S_n$, suppose that there
exist $i_1$, $i_2$, $\ldots$, $i_6 \in [n]$ such that
\begin{eqnarray*} \gamma_1(i_1)&=& i_2,~ \gamma_1(i_2)=i_1, ~\gamma_1(i_3)=i_4,~ \gamma_1(i_4)=i_3, ~\gamma_1(i_5)=i_5,~ \gamma_1(i_6)=i_6, \\
\gamma_2(i_1)&=& i_2,~ \gamma_2(i_2)=i_1, ~\gamma_2(i_3)=i_3,~ \gamma_2(i_4)=i_4, ~\gamma_2(i_5)=i_5,~ \gamma_1(i_6)=i_6, \end{eqnarray*}
and for all $j \in [n] \setminus \{i_1$, $i_2$, $\ldots$, $i_6 \}$, $\gamma_1(j)=\gamma_2 (j)$. Then, one can check that
 all equipartitions  fixed by $\gamma_1$ are also fixed by $\gamma_2$. This
implies that  $\alpha$ is a  transposition or a product of $m$
disjoint transpositions. Note that the number of equipartitions
fixed by a transposition is ${2m-2 \choose m-2}$ and the number of
equipartitions fixed by a product of $m$ disjoint transpositions
is $2^{m-1}$ if $m$ is odd; $2^{m-1} + {m \choose m/2}$ if $m$ is
even. Since for $m \ge 4$, ${2m-2 \choose m-2} > 2^{m-1} + {m
\choose m/2}$, $\alpha$ is a transposition. Therefore for any
non-identity permutation $\beta \in S_n$, the number of
equipartitions fixed by $\beta$ is at most ${2m-2 \choose m-2}$.

\begin{lem} \label{random-3-coloring}
For $n=2m$ with $m \ge 4$, let $\Phi$ be the set of equipartitions
of $[n]$. Then there is a 3-coloring $c: \Phi \rightarrow \{B,R,Y
\}$ of $\Phi$ such that only identity permutation in $S_n$
preserves all colors under the induced action of $S_n$ on $\Phi$.
\begin{proof}
Give a random coloring on $\Phi$ with three colors $\{B,R,Y \}$.
For any non-identity permutation $\beta$ of $[2m]$, let
$A_{\beta}$ be the event that $\beta$ preserves colors of all
equipartitions of $[2m]$. Note that the number of equipartitions
fixed by $\beta$ is at most ${2m-2 \choose m-2}$. Namely,  the
number of equipartitions which are not fixed by $\beta$ is at
least
$$|\Phi|
- {2m-2 \choose m-2} = \frac{1}{2}{2m \choose m} - {2m-2 \choose
m-2} = \frac{m}{m-1}{2m-2 \choose m-2}.$$
 For any orbit $O$ of $\beta$ whose
size is $t$ with $t \ge 2$ under the action of $S_{2m}$ on $\Phi$,
the probability that $\beta$ preserves colors of all
equipartitions in $O$ is $3^{-t+1}$, which is less than $3^{-\frac{t}{2}}$. Hence, we have
$$Pr(A_{\beta}) \le 3^{-
\frac{m}{2(m-1)}{2m-2 \choose m-2}}.$$ Therefore,
$$Pr(\bigcup_{\beta \in S_n \setminus \{ id \}} A_{\beta})) \le \sum_{\beta \in S_n \setminus \{ id \}} Pr(A_{\beta})
< ((2m)!-1) 3^{- \frac{m}{2(m-1)}{2m-2 \choose m-2}} \le (2m)!
3^{- \frac{m}{2(m-1)}{2m-2 \choose m-2}},$$ where $id$ is the
identity permutation of $[n]$. For $n=8$, the number $n! 3^{-
\frac{m}{2(m-1)}{2m-2 \choose m-2}}$ is $\dfrac{4480}{3^8}$, and it
is less than 1. Furthermore, for any $m \ge 4$, $$\frac{(2m+2)!
3^{- \frac{m+1}{2m}{2m \choose m-1}}}{(2m)! 3^{-
\frac{m}{2(m-1)}{2m-2 \choose m-2}}} = \frac{(2m+2)(2m+1)}{ 3^{
\frac{3m-2}{2(m-1)}{2m-2 \choose m-2}}} < \frac{(2m+2)(2m+1)}{ 3^{
\frac{3}{2}{2m-2 \choose m-2}}}<1$$ because
\begin{eqnarray*} 3^{ \frac{3}{2}{2m-2 \choose m-2}} &=& 3^{ \frac{3(2m-2)(2m-3)\cdots (m+1)}{2\cdot (m-2)!}}
= 3^{ \frac{(m-1)(2m-3)(2m-5)(2m-6)\cdots (m+1)}{(m-3)(m-4)\cdots
5\cdot4}} \\ &>& 3^{ (m-1)(2m-3)(2m-5)} > (2m+2)(2m+1)
\end{eqnarray*} for any $m \ge 4$. This implies that for any $m \ge
4$,
$$Pr(\bigcup_{\beta \in S_n \setminus \{ id \}} A_{\beta})) < 1$$
and hence there exists a $3$-coloring $c: \Phi \rightarrow \{B,R,Y
\}$ of $\Phi$ such that the identity permutation in $S_{2m}$ is the only
color-preserving permutation under the induced action of $S_{2m}$
on $\Phi$.
\end{proof}
\end{lem}

\section{A Proof of the main theorem} \label{proof}

In this section, we prove Theorem \ref{main} which is the main
result in this paper. For a graph $G$ and for a $v \in V(G)$, let
$N(v)$ be the set of all vertices adjacent to $v$. For any
positive integer $i$, let $D_i$ be the set of all vertices whose
degrees are $i$.

Let $J=J(n,k)_I$, where $1 \le k \le \frac{n}{2}$ and $\emptyset
\varsubsetneq I \varsubsetneq \{1,2,\ldots,k\}$. If $k=1$ then the
graph $J$ is a complete graph $K_n$, and hence its distinguishing
number is $n$. Assume that $k > 1$.

\begin{figure}
$$
\begin{pspicture}[shift=-.9](-7.3,-5.5)(7.3,7.5)
\pscircle(0,0){5}
\psline(5;0)(5;90)(5;180)(5;270)(5;0)
\psline(5;30)(5;120)(5;210)(5;300)(5;30)
\psline(5;60)(5;150)(5;240)(5;330)(5;60)
\psline(0;90)(7;45)
\psline(0;90)(5;-60)
\psline(0;90)(5;180)
\psline(0;90)(5;150)
\psline(5;90)(5;-90)
\psline(3;90)(5;60)
\psline(3;90)(5;30)
\psline(3;90)(5;-30)
\psline(3;90)(5;150)
\psline(3;90)(5;180)
\psline(2;165)(5;90)
\psline(2;165)(5;150)
\psline(2;165)(5;180)
\psline(2;165)(5;210)
\psline(2;165)(5;270)
\psline(2;165)(7;45)
\psline(7;45)(5;0)
\psline(7;45)(5;30)
\psline(7;45)(5;60)
\psline(7;135)(5;120)
\psline(7;135)(5;150)
\psline(7;135)(5;180)
\pccurve[angleA=45,angleB=90](7;135)(5;30)
\pccurve[angleA=-135,angleB=180](7;135)(5;240)
\pccurve[angleA=135,angleB=90](7;45)(5;120)
\pccurve[angleA=-45,angleB=0](7;45)(5;-30)
\pscircle[fillstyle=solid,fillcolor=darkgray](0;90){.35}
\pscircle[fillstyle=solid,fillcolor=darkgray](5;90){.35}
\pscircle[fillstyle=solid,fillcolor=darkgray](5;60){.35}
\pscircle[fillstyle=solid,fillcolor=darkgray](5;30){.35}
\pscircle[fillstyle=solid,fillcolor=darkgray](5;0){.35}
\pscircle[fillstyle=solid,fillcolor=darkgray](5;-30){.35}
\pscircle[fillstyle=solid,fillcolor=darkgray](5;-60){.35}
\pscircle[fillstyle=solid,fillcolor=darkgray](5;-90){.35}
\pscircle[fillstyle=solid,fillcolor=darkgray](5;-120){.35}
\pscircle[fillstyle=solid,fillcolor=darkgray](5;-150){.35}
\pscircle[fillstyle=solid,fillcolor=darkgray](5;180){.35}
\pscircle[fillstyle=solid,fillcolor=darkgray](5;150){.35}
\pscircle[fillstyle=solid,fillcolor=darkgray](5;120){.35}
\pscircle[fillstyle=solid,fillcolor=darkgray](3;90){.37}
\pscircle[fillstyle=solid,fillcolor=darkgray](2;165){.37}
\pscircle[fillstyle=solid,fillcolor=darkgray](7;135){.37}
\pscircle[fillstyle=solid,fillcolor=darkgray](7;45){.37}
\rput(5;90){$V_1$}
\rput(5;60){$V_2$}
\rput(5;30){$V_3$}
\rput(5;0){$V_4$}
\rput(5;-30){$V_5$}
\rput(5;-60){$V_6$}
\rput(5;-90){$V_7$}
\rput(5;-120){$V_8$}
\rput(5;-150){$V_9$}
\rput(5;180){$V_{10}$}
\rput(5;150){$V_{11}$}
\rput(5;120){$V_{12}$}
\rput(3;90){$V_{2}^{\tau_2}$}
\rput(2;165){$V_{10}^{\tau_2}$}
\rput(7;45){$V_{3}^{\tau_2}$}
\rput(7;135){$V_{11}^{\tau_2}$}
\rput(1.25;-4){$\{1,3,5,7\}$}
\end{pspicture}
$$
\caption{.}  \label{fig1}
\end{figure}

\medskip

\noindent Case 1: $(n,k)\neq (5,2)$ and $2 \le k < \frac{n-1}{2}$,
and $J \neq J(12,4)_I$ with $I= \{1,3 \}$ or  $\{2,4 \}$.

In this case, $\Aut(J) = \Aut(K(n,k))$ as an acting group on the
vertex set, and hence $\Dist(J) = \Dist(K(n,k))=2$ by Proposition
\ref{dist-Kneser} and Lemma \ref{same-dist}.

\medskip

\noindent Case 2: $J = J(12,4)_I$ with $I= \{1,3 \}$ or  $\{2,4
\}$.

Assume that $J = J(12$, $4)_I$ with $I= \{1$, $3 \}$. Let $S_1 = \{
V_{j}$, $V^{\tau_1}_{j}$, $V^{\tau_2}_{j} \ | \ j=1$, $2$, $\ldots$, $12 \}$ $
\cup$ $\{ X_1 =\{1$, $3$, $5$, $7 \} \}$ and let $H_1$ be a subgraph of $J$
induced by $S_1$ as illustrated in Figure~\ref{fig1}. Then $S_1$ is a determining set
by Lemma \ref{12-4--det}. Let $\psi_1$ be an automorphism of
$H_1$.  Note that the order of $H_1$ is $17$ and
\begin{eqnarray*}
D_5&=& \{V_4 , V_5, V_8, V_9 \},~~D_6 = \{V_1, V_6, V_7, V_{12},
V^{\tau_2}_{11}, X_1 \} \\
D_7&=& \{ V_3, V^{\tau_1}_{2}, V^{\tau_1}_{10} \},~~D_8 = \{V_2,
V_{10},  V^{\tau_2}_{3} \},~~D_9 =\{ V_{11} \}
\end{eqnarray*}
Since $D_9 = \{ V_{11} \}$, $\psi_1(V_{11})=V_{11}$.  The fact
$N(V_{11}) \cap D_5 = \{ V_8 \}$  implies that $\psi_1(V_8)=V_8$.
Since $N(V_8) \cap D_5 = \{ V_5, V_9 \}$,  $\psi_1$ fixes $V_4$.
By the fact that $N(V_4) \cap N(V_8) \cap D_5 = \{ V_5 \}$,
$\psi_1$ fixes both $V_5$ and $V_9$.  Note that $N(V_8) \cap
N(V_{11}) = \{V^{\tau_2}_{11} \}$. So $\psi_1$ fixes
$V^{\tau_2}_{11}$.  The fact $N(V_4) \cap N(V_8) \cap D_6 = \{V_7
\}$ and $N(V_5) \cap N(V_9) \cap D_6 = \{V_6 \}$ implies that
$\psi_1$ fixes both $V_6$ and $V_7$. Furthermore, $\psi_1$ fixes
$X_1=\{1$, $3$, $5$, $7 \}$ because $N(V_6 ) \cap N(V_7) = \{ X_1 \}$.
Since $N(V_7)$ $\cap$ $N(V_{11})$ $\cap$ $D_7$ $= \{ V^{\tau_1}_{10} \}$ and
$N(V_7)$ $\cap$ $N(V_{11})$ $\cap$ $D_8 =$ $\{ V_{10} \}$, $\psi_1$ also
fixes both $V_{10}$ and $V^{\tau_1}_{10}$. Note that $N(V_4) \cap
N(V_6) \cap D_7$ $= \{ V_3 \}$ and $N(V_4)$ $\cap N(V_6)$ $\cap D_8$ $= \{
V^{\tau_2}_{3} \}$. This implies that $\psi_1(V_3)$ $=V_3$ and
$\psi_1(V^{\tau_2}_{3})$ $=V^{\tau_2}_{3}$.  By the fact $N(V_5) \cap
N(V_{11})$ $= \{ V^{\tau_1}_{2} \}$ and $N(V_9) \cap
N(V^{\tau_2}_{11})$ $\cap D_6$ $= \{ V^{12} \}$, $\psi_1$ fixes both
$V^{\tau_1}_{2}$ and $V_{12}$. Up to now, we know that $\psi_1$
fixes all vertices in $V(H)$ except $V_1$ and $V_2$. Since the
degree of $V_1$ is 6 and that of $V_2$ is 8, $\psi_1$ also fixes
both $V_1$ and $V_2$. Therefore, $\psi_1$ is the identity, which
implies that $H$ is an asymmetric graph. By Corollary
\ref{asymm-dist}, we have $\Dist(J) = 2$.

For any $J_1 = J(12,4)_I$ with  $I= \{2,4 \}$,
$\Aut(J_1)=\Aut(J)$. Therefore, $\Dist(J_1) = \Dist(J)=2$.

\medskip

\noindent Case 3: $J = J(5,2)_I$ with $I= \{1 \}$ or  $\{2 \}$.

Since $J(5,2)_{\{ 2 \}}$ is the Kneser graph $K(5,2)$ and
$J(5,2)_{\{ 1 \}}$ is its complement, we have $\Dist(J)=3$ by Proposition
\ref{dist-Kneser}.

\medskip

\noindent Case 4: $k=\frac{n-1}{2}$ and $I \neq k+1-I$.

If $n=5$ then $I= \{1 \}$ or  $\{2 \}$, which means that $I =
k+1-I$. Hence we have $n \ge 7$.  In this case, since $\Aut(J) =
\Aut(K(n,\frac{n-1}{2}))$ as an acting group on the vertex set,
one can find that $\Dist(J) = \Dist(K(n,k))=2$ by Proposition~\ref{dist-Kneser} and Lemma~\ref{same-dist}.

\begin{figure}
$$
\begin{pspicture}[shift=-.9](-.5,-1)(12,2.5)
\psline(0,1.5)(5,1.5)
\psline(6,1.5)(11,1.5)
\psline(0,0)(0,1.5)
\psline(1.5,0)(1.5,1.5)
\pccurve[angleA=-35,angleB=-145](0,1.5)(11,1.5)
\pccurve[angleA=-25,angleB=-155](0,1.5)(9.5,1.5)
\pccurve[angleA=-25,angleB=-155](1.5,1.5)(11,1.5)
\pscircle[fillstyle=solid,fillcolor=darkgray](0,1.5){.2}
\pscircle[fillstyle=solid,fillcolor=darkgray](1.5,1.5){.2}
\pscircle[fillstyle=solid,fillcolor=darkgray](3,1.5){.2}
\pscircle[fillstyle=solid,fillcolor=darkgray](4.5,1.5){.2}
\pscircle[fillstyle=solid,fillcolor=darkgray](6.5,1.5){.2}
\pscircle[fillstyle=solid,fillcolor=darkgray](8,1.5){.2}
\pscircle[fillstyle=solid,fillcolor=darkgray](9.5,1.5){.2}
\pscircle[fillstyle=solid,fillcolor=darkgray](11,1.5){.2}
\pscircle[fillstyle=solid,fillcolor=darkgray](0,0){.2}
\pscircle[fillstyle=solid,fillcolor=darkgray](1.5,0){.2}
\rput(0,2){$V_1$}
\rput(1.5,2){$V_2$}
\rput(3,2){$V_3$}
\rput(4.5,2){$V_4$}
\rput(6.5,2){$V_{m-1}$}
\rput(8,2){$V_m$}
\rput(9.5,2){$V_{m+1}$}
\rput(11,2){$V_{m+2}$}
\rput(5.5,1.5){$\cdots$}
\rput(0,-.5){$X_2$}
\rput(1.5,-.5){$Y_2$}
\end{pspicture}
$$
\caption{.}  \label{fig2}
\end{figure}

\medskip

\noindent Case 5: $k=\frac{n-1}{2}$ and $I = k+1-I$.

Let $J= J(2m+1,m)_{\{1,m\}}$ with $m \ge 3$. Let
$$S_2 = \{ V_1, V_2, \ldots, V_{m+2} \} \cup \{ X_2, Y_2 \},$$ where
$X_2=\{1$, $2$, $\ldots$, $m-2$, $m$, $m+2 \} $ and $Y_2 =\{2$, $3$, $\ldots$, $m-1$,
$m+1$, $m+3 \}$; and let $H_2$ be a subgraph of $J$ induced by
$S_2$ as depicted in Figure~\ref{fig2}. Then, $S_2$ is a determining set of $J$ by
Lemma~\ref{n,k-det}(1). Let $\psi_2$ be an automorphism of $J$.
Note that
$$D_{1} = \{ X_2 , Y_2 \},~~D_2 =\{ V_3, V_4, \ldots, V_m \},
~~D_3 = \{ V_{m+1}, V_{m+2}\},~~D_4 = \{ V_1, V_2 \}.$$ Since
$V_1$ is the only vertex adjacent to all vertices in $D_3$,
$\psi_2(V_1) = V_1$. This implies that $\psi_2$ fixes $V_2, X_2$
and $Y_2$ because $D_4 = \{ V_1$, $V_2 \}$; and $X_2$ and $Y_2$ are
only adjacent to $V_1$ and $V_2$, respectively. By the fact that
$V_{m+2}$ is the only vertex adjacent to all vertices in $D_4$,
$\psi_2(V_{m+2}) = V_{m+2}$. This implies that $\psi_2$ also fixes
$V_{m+1}$ because $D_3 = \{ V_{m+1}$, $V_{m+2}\}$. Since $N(V_2)
\cap D_2 = \{ V_3 \}$, $\psi_2 (V_3) =V_3$. By a similar way, one
can show that $\psi_2$ fixes all vertices in $H_2$, and hence
$\psi_2$ is the identity. This means that $H_2$ is asymmetric. By
Corollary \ref{asymm-dist}, we have $\Dist(J) = 2$.

For any $J_2 = J(2m+1,m)_I$ satisfying  $I = k+1-I$,
$\Aut(J_2)=\Aut(J)$. Therefore $\Dist(J_2) = \Dist(J)=2$.

\medskip

\noindent Case 6: $k=\frac{n}{2}$ and $I$ is neither $\{k\}$ nor
$\{1,2,\ldots,k-1\}$, and $I' \neq I''$.

Let $J= J(2m,m)_{\{1\}}$ with $m \ge 3$. Let
$$S_3 = \{ V_1, V_2, \ldots, V_{2m} \} \cup \{ X_3, Y_3, Z_3 \},$$ where
$X_3=\{1$, $2$, $\ldots$, $m-2$, $m$, $m+2 \} $, $Y_3 =\{2$, $3$, $\ldots$, $m-1$, $m+1$,
$m+3 \}$ and $Z_3 =\{4$, $5$, $\ldots$, $m+1$, $m+3$, $m+5 \}$; and let $H_3$
be a subgraph of $J$ induced by $S_3$ as shown in Figure~\ref{fig3}. Then, $S_3$
is a determining set of $J$ by Lemma \ref{n,k-det}(2). Let
$\psi_3$ be an automorphism of $J$.  Note that
$$D_{1} = \{ X_3 , Y_3, Z_3 \},~~D_3 =\{ V_1, V_2,  V_4 \}
~~\mbox{and}~~D_2 = V(H_3) - (D_1 \cup D_3).$$  By a similar way
with Cases 2 and 5, one can show that $\psi_3$ fixes all vertices
in $H_3$, and hence $\psi_3$ is the identity. So $H_3$ is an
asymmetric graph. By Corollary \ref{asymm-dist}, we have $\Dist(J) = 2$.

For any $J_3 = J(2m,m)_I$ satisfying  $I \neq \{k\}$ nor
$\{1,2,\ldots,k-1\}$, and $I' \neq I''$, $\Aut(J_3)=\Aut(J)$.
Therefore, we find $\Dist(J_3) = \Dist(J)=2$.

\begin{figure}
$$
\begin{pspicture}[shift=-.9](-.5,-1)(12,2.5)
\psline(0,1.5)(5,1.5)
\psline(6,1.5)(11,1.5)
\psline(0,0)(0,1.5)
\psline(1.5,0)(1.5,1.5)
\psline(4.5,0)(4.5,1.5)
\pccurve[angleA=-20,angleB=-160](0,1.5)(11,1.5)
\pscircle[fillstyle=solid,fillcolor=darkgray](0,1.5){.2}
\pscircle[fillstyle=solid,fillcolor=darkgray](1.5,1.5){.2}
\pscircle[fillstyle=solid,fillcolor=darkgray](3,1.5){.2}
\pscircle[fillstyle=solid,fillcolor=darkgray](4.5,1.5){.2}
\pscircle[fillstyle=solid,fillcolor=darkgray](6.5,1.5){.2}
\pscircle[fillstyle=solid,fillcolor=darkgray](8,1.5){.2}
\pscircle[fillstyle=solid,fillcolor=darkgray](9.5,1.5){.2}
\pscircle[fillstyle=solid,fillcolor=darkgray](11,1.5){.2}
\pscircle[fillstyle=solid,fillcolor=darkgray](0,0){.2}
\pscircle[fillstyle=solid,fillcolor=darkgray](1.5,0){.2}
\pscircle[fillstyle=solid,fillcolor=darkgray](4.5,0){.2}
\rput(0,2){$V_1$}
\rput(1.5,2){$V_2$}
\rput(3,2){$V_3$}
\rput(4.5,2){$V_4$}
\rput(6.5,2){$V_{2m-3}$}
\rput(8,2){$V_{2m-2}$}
\rput(9.5,2){$V_{2m-1}$}
\rput(11,2){$V_{2m}$}
\rput(5.5,1.5){$\cdots$}
\rput(0,-.5){$X_3$}
\rput(1.5,-.5){$Y_3$}
\rput(4.5,-.5){$Z_3$}
\end{pspicture}
$$
\caption{.}  \label{fig3}
\end{figure}

\medskip

\noindent Case 7: $k=\frac{n}{2}$, $I' = I''$ and $I$ is neither $\{k\}$
nor $\{1$, $2$, $\ldots$, $k-1\}$.

Note that for $n \le 6$, this case can not occur. Hence assume
that $n \ge 8$.  Let $J= J(2m,m)_{\{1, m-1 \}}$ with $m \ge 4$. By
Lemma~\ref{n=2m-m-Dist>2}, we have $\Dist(J) \ge 3$. Let $f: V(J)
\rightarrow \{1$, $2$, $3 \}$ be a random coloring satisfying that for
any $u \in V(J)$, $f(u)$ and $f(\bar{u})$ are distinct. Note that
for any $u \in V(J)$,
$$Pr(\{f(u), f(\bar{u})\} = \{1,2 \}) = Pr(\{f(u), f(\bar{u})\} = \{1,3
\})=Pr(\{f(u), f(\bar{u})\} = \{2,3 \})=\frac{1}{3}.$$

Let $\Phi$ be the set of equipartitions of $[n]$. Using a random
coloring $f: V(J) \rightarrow \{1$, $2$, $3 \}$, we can define the
coloring $\tilde{f}: \Phi \rightarrow \{B$, $R$, $Y \}$ as follows:
for any $\{u$, $\bar{u} \} \in \Phi$,  $\tilde{f}(\{u$, $\bar{u} \})
= B$ if $\{f(u)$, $f(\bar{u})\}$ $= \{1,2 \}$; $\tilde{f}(\{u$, $
\bar{u} \}) = R$ if $\{f(u)$, $f(\bar{u})\} = \{1$, $3 \}$; and
$\tilde{f}(\{u$, $\bar{u} \}) = Y$ if $\{f(u)$, $f(\bar{u})\} = \{2$, $3
\}$. Then we can consider $\tilde{f}$ as a random 3-coloring of
$\Phi$. By Lemma~\ref{random-3-coloring}, there is a random
$3$-coloring $f: V(J) \rightarrow \{1$, $2$, $3 \}$ such that only
identity permutation in $S_n$ preserves all colors of
equipartitions in its corresponding random $3$-coloring $\tilde{f}:
\Phi \rightarrow \{B$, $R$, $Y \}$ under the induced action of $S_n$
on $\Phi$. This implies that there exists a $3$-distinguishing
coloring $f: V(J) \rightarrow \{1$, $2$, $3 \}$. Therefore, $\Dist(J)
\le 3$, and hence $\Dist(J)=3$.

For any $J_4 = J(2m,m)_I$ satisfying  $I' = I''$ and $I$ is neither
$\{k\}$ nor $\{1$, $2$, $\ldots$, $k-1\}$, $\Aut(J_4)=\Aut(J)$. Therefore,
$\Dist(J_4) = \Dist(J)=3$.

\medskip

\noindent Case 8: $k=\frac{n}{2}$ and $I = \{k\}$ or
$\{1,2,\ldots,k-1\}$.

Let $J= J(2m,m)_{\{m \}}$. Then $J$ is composed of $ \frac{{2m
\choose m}}{2}$ components which are isomorphic to $K_2$. Note
that a coloring $f: V(J) \rightarrow \{1$, $2$, $\ldots$, $r\}$ is an
$r$-distinguishing if and only if for any vertex $u \in V(J)$,x
$f(u)$ and $f(\bar{u})$ are distinct and for any two vertex $v$, $w
\in V(J)$ contained to different components, $\{ f(v)$, $f(\bar{v})
\} \neq \{ f(w)$, $f(\bar{w}) \}$. Hence $\Dist(J)$ is the smallest
integer $r$ such that ${r \choose 2} \ge \frac{{2m
\choose m}}{2}$. Therefore, $\Dist(J) = \lceil \frac{1+
\sqrt{1+4{2m \choose m}}}{2} \rceil$.

For any $J_5 = J(2m,m)_I$ with $I = \{1$, $2$, $\ldots$, $m-1\}$, $J_5$ is
the complement of $J$. Hence, $\Dist(J_5) = \Dist(J)=\lceil
\frac{1+ \sqrt{1+4{2m \choose m}}}{2} \rceil$.

\section*{Acknowledgments}
The author would like to thank ....
The \TeX\, macro package
PSTricks~\cite{PSTricks} was essential for typesetting the equations
and figures.

 \end{document}